\documentclass[10pt]{article}
\usepackage{amssymb,amsbsy,amsmath,amsfonts,amssymb,amscd, mathrsfs}

\usepackage[english]{babel}
\usepackage[T1]{fontenc}
\usepackage{indentfirst}

\makeatletter
\@addtoreset{equation}{section}
\makeatother

\newtheorem{statement}{}[section]
\newtheorem{theoreme}[statement]{Theorem}
\newtheorem{lemme}[statement]{Lemma}

\newtheorem{proposition}[statement]{Proposition}
\newtheorem{definition}[statement]{Definition}
\newtheorem{corollaire}[statement]{Corollary}

\newcommand\C{\mathbb C}
\newcommand\N{\mathbb N}
\newcommand\R{\mathbb R}
\newcommand\T{\mathbb T}
\newcommand\D{\mathbb D}
\newcommand\Z{\mathbb Z}
\newcommand\e{{\rm e}}
\renewcommand\P{\mathbb P}
\newcommand\eps{\varepsilon}

\let\hat=\widehat
\newcommand\ind{{\rm 1\kern-.30em I}}
\newcommand\qed{\hfill $\square$}

\let\phi=\varphi

\title{\bf Some translation-invariant Banach function spaces which contain $c_0$}
\author{\it Pascal Lef\`evre, Daniel Li,\\ \it Herv\'e Queff\'elec, Luis Rodr{\'\i}guez-Piazza}

\date{}

\begin{document}

\maketitle

\noindent{\bf Abstract.} \emph{We produce several situations where some natural subspaces of classical Banach 
spaces of functions over a compact abelian group contain the space $c_0$.} 
\medskip

\medskip



\section{Introduction}

Let $G$ be a compact abelian group and $\Gamma=\hat G$ its dual group. It is a familiar theme in Harmonic 
Analysis to compare the ``thinness'' properties of a subset $\Lambda\subseteq \Gamma$ with the Banach space
properties of the space $X_\Lambda$, where $X$ is a Banach space of Haar-integrable functions on $G$ and 
$X_\Lambda$ is the subspace of $X$ consisting of the $f\in X$ whose spectrum lies in $\Lambda$:
$\hat f (\gamma) = 0$ if $\gamma\notin \Lambda$. We refer to Kwapie\'n-Pe{\l}czy\'nski's classical paper 
\cite{K-P} for such investigations.
\par
It is known that, denoting by $\Psi_2$ the Orlicz function $\e^{x^2}-1$:\par
\medskip

$(1)$ If $L^{\Psi_2}_\Lambda = L^2_\Lambda$, then $\Lambda$ is a Sidon set (Pisier \cite{Pisier}, 
Th\'eor\`eme 6.2);\par
\smallskip

$(2)$ If ${\cal C}_\Lambda$ has a finite cotype, then $\Lambda$ is a Sidon set (Bourgain-Milman 
\cite{Bourgain-Milman}).
\par\medskip

\noindent
Recall that $\Lambda$ is a \emph{Sidon set} if every continuous function on $G$ with spectrum in $\Lambda$ 
has an absolutely convergent Fourier series.
\par\smallskip

In a previous paper, we proved, among other facts, the following extension of $(1)$ (\cite{LLQR}, Theorem 2.3):
\par\medskip

$(1')$ If $L^{\Psi_2}_\Lambda$ has cotype $2$, then $\Lambda$ is a Sidon set;\par
\smallskip

We also showed the following variant of $(2)$ (\cite{LLQR}, Theorem 1.2):\par
\medskip

$(2')$ If $U_\Lambda$ has a finite cotype, then $\Lambda$ is a Sidon set,\par
\smallskip

\noindent
where $U=U(\T)$ is the space of the continuous functions on the circle group $\T$ whose Fourier series converges 
uniformly on $\T$.\par
\medskip

In this work, we study what are the implications on $\Lambda$ of the fact that some Banach space 
$X_\Lambda$ contains, or not, the space $c_0$. In particular, we shall extend $(1')$ and $(2')$.
\par
\medskip

The paper is organized as follows. In Section 2, we show that if $\psi$ is an Orlicz function which violates the 
$\Delta_2$-condition, in a strong sense: $\displaystyle \lim_{x\to+\infty} \frac{\psi(2x)}{\psi(x)}=+\infty$
(which is the case of $\Psi_2$), and if $X_0$ is a linear subspace of $L^\infty$ on which the norms 
$\Vert\ \Vert_2$ and $\Vert\ \Vert_\psi$ are not equivalent, then the closure $X$ of $X_0$ in $L^\psi$ 
contains $c_0$. It follows that if $\Lambda$ is not a Sidon set, then $L^{\Psi_2}_\Lambda$ contains $c_0$, 
and {\it a fortiori} that if $L^{\Psi_2}_\Lambda$ has a finite cotype, then $\Lambda$ is a Sidon set, which 
generalizes $(1')$.\par
In Section 3, we extend $(2')$ by showing that: If $\Lambda$ is not a set of uniform convergence ({\it i.e.} if
$U_\Lambda\not={\cal C}_\Lambda$), then $U_\Lambda$ does contain $c_0$. In particular, if $U_\Lambda$ 
has a finite cotype, then $U_\Lambda={\cal C}_\Lambda$, so ${\cal C}_\Lambda$ has a finite cotype and
therefore, in view of $(2)$, $\Lambda$ is a Sidon set. This explains why the proof of $(2')$ in \cite{LLQR} 
mimicked Bourgain and Milman's.\par
In Section 4, we use the notion of invariant mean in $L^\infty(G)$. We say that $\Lambda$ is a 
\emph{Lust-Piquard set} if, for every function $f\in L^\infty_\Lambda$, the product $\gamma f$ of $f$ with every 
character $\gamma\in \Gamma$ has a unique invariant mean. Of course, if every $f\in L^\infty_\Lambda$ is 
continuous ({\it i.e.} $\Lambda$ is a Rosenthal set), then $\Lambda$ is a Lust-Piquard set. F. Lust-Piquard 
(\cite{LP3}) showed that there are Lust-Piquard sets which are not Rosenthal sets, and, more precisely, that 
$\Lambda={\P}\cap (5\Z+2)$, where ${\P}$ is the set of the prime numbers, is a Lust-Piquard set such that 
${\cal C}_\Lambda$ contains $c_0$ (if $\Lambda$ is a Rosenthal set, ${\cal C}_\Lambda$ cannot
contain $c_0$). We construct here another kind of ``big'' Lust-Piquard set $\Lambda$, namely a Hilbert set. 
Then ${\cal C}_\Lambda$ contains $c_0$ by a result of the second-named author (\cite{Li2}, Theorem 2).\par
In Section 5, we investigate under which conditions the space ${\cal C}_\Lambda$ is complemented in 
$L^\infty_\Lambda$. We conjecture that this happens only if ${\cal C}_\Lambda=L^\infty_\Lambda$, {\it i.e.}
$\Lambda$ is a Rosenthal set. We are only able to show that, under that condition of complementation, 
${\cal C}_\Lambda$ does not contain $c_0$, and, moreover, every $f\in L^\infty_\Lambda$ which is 
Riemann-integrable is actually in ${\cal C}_\Lambda$.
\medskip

\noindent{\bf Notation.} Throughout this paper, $G$ is a compact abelian group, and $\Gamma=\hat G$ is its 
(discrete) dual group. The Haar measure of $G$ is denoted by $m$, and integration with respect to $m$ by $dt$ or 
$dx$. We shall write the group structure of $\Gamma$ additively, so that, for $\gamma\in \Gamma$, the 
character $(-\gamma)\in \Gamma$ is the function $\overline\gamma\in {\cal C}(G)$. When $G$ is the circle 
group $\T=\R/2\pi\Z$, we identify, as usual, the character $e_n\colon t\mapsto \e^{int}$ with the
integer $n\in \Z$, and so the dual group $\Gamma$ to $\Z$; the Haar measure is then $dt/2\pi$.\par
For $f\in L^1(G)$, the Fourier coefficient of $f$ at $\gamma\in \Gamma$ is 
$\hat f(\gamma)=\int_G f(t)\,\overline{\gamma(t)}\,dt$. If $X$ is a linear function subspace of 
$L^1(G)$, we denote by $X_\Lambda$ the subspace of those $f\in X$ for which the Fourier coefficients vanish 
outside of $\Lambda$.\par
When we say that a Banach space $X$ contains a Banach space $Y$, we mean that $X$ contains a (closed) 
subspace isomorphic to $Y$.

\medskip
\noindent{\bf Aknowledgements.} This work was partly supported by a Picasso project (EGIDE-MCYT) between 
the french and spanish governments.\par
We thank the referee for a very careful reading of this paper and for many suggestions to improve the writing.

\section{Subspaces of Orlicz spaces}

Let \hfill $\psi$ \hfill be \hfill an \hfill Orlicz \hfill function, \hfill that \hfill is, \hfill an 
\hfill increasing \hfill convex \hfill function\\ 
$\psi\colon [0,+\infty\,[ \to [0,+\infty\,[$ 
such that $\psi(0)=0$ and $\psi(+\infty)=+\infty$. We shall assume that $\psi$ violates the $\Delta_2$-condition, 
in the following strong sense:
\begin{equation}
\lim_{x\to+\infty} \frac{\psi(2x)}{\psi(x)} = +\infty\,.\tag*{($\ast$)}
\end{equation}

Let $(\Omega, {\cal A}, \P)$ be a probability space. The Orlicz space $L^\psi(\Omega)$ is the space of all the 
(equivalence classes of) measurable functions $f\colon \Omega\to\C$ for which there is a constant $C\geq 0$ such
that
$$\int_\Omega \psi\Big(\frac{\vert f(t)\vert}{C}\Big)\,d\P (t)\leq 1$$
and then $\Vert f\Vert_\psi$ is the least possible constant $C$.\par
Observe that $(\ast)$ implies that there exists $a>0$ such that $\psi(2t)\geq 4\,\psi(t)$ for every $t\geq a$. 
Hence, for all $n\geq 0$, one has $\psi(2^n a)\geq 4^n\psi(a)$. It follows that, for
$2^n a\leq x < 2^{n+1} a$, we have
$$\psi(x)\geq \psi(2^n a)\geq 4^n\psi(a)\geq \Big(\frac{x}{2a}\Big)^2\psi(a) =C\,x^2\,.$$
Hence $\psi (x) \geq C\,x^2$ for every $x\geq a$, and so the norm $\Vert\ \Vert_\psi$ is stronger than the norm 
of $L^2$.
\medskip

\begin{theoreme}\label{theo 1}
Suppose that $\psi$ is an Orlicz function as above. Let $X_0$ be a linear subspace of $L^\infty(\Omega)$ on 
which the norms $\Vert\ \Vert_2$ and $\Vert\ \Vert_\psi$ are not equivalent. Then there exists in $X_0$ a 
sequence which is equivalent, in the closure $X$ of $X_0$ for the norm $\Vert\ \Vert_\psi$, to the canonical 
basis of $c_0$.
\end{theoreme}

\noindent{\bf Proof.} We first remark that, thanks to $(\ast)$, we can choose, for each $n\geq 1$, a positive number 
$x_n$ such that
$$\hskip 1cm\psi\Big(\frac{x}{2}\Big)\leq \frac{1\,}{\,2^n}\,\psi(x)\,,
\hskip 5mm\forall x\geq x_n\,.$$
Since $\psi$ increases, we have for every $x\geq 0$:
$$\hskip 3mm
\psi\Big(\frac{x}{2}\Big)\leq \frac{1\,}{\,2^n}\,\psi(x)+\psi(x_n)\,.$$\par
Next, $\psi$ is continuous since it is convex. Hence there exists $a>0$ such that $\psi(a)=1$. Then, since $\psi$ is 
increasing, we have, for every $f\in L^\infty(\Omega)$:
$$\int_\Omega
\psi\Big(a\, \frac{\vert f\vert\hskip 3mm}{\Vert f\Vert_\infty}\Big) \,d\P \leq 1\,,$$
and so $\Vert f\Vert_\psi\leq (1/a)\,\Vert f\Vert_\infty$.\par
\smallskip
Now, let $\alpha_n$, $n\geq 1$, be positive numbers less than $a/2$ such that $\sum_{n\geq 1} \alpha_n < a$.
We shall construct inductively a sequence of functions $f_n\in X_0$, with \! $\Vert f_n\Vert_\psi=1$, and a 
sequence of positive numbers $\beta_n\leq 1/2^n$ such that:
\begin{itemize}
\item[(i)] $\P (\{\vert f_n\vert>\alpha_n\})\leq \beta_n$, for every $n\geq 1$;\par
\item[(ii)] if we set $M_1=1$ and, for $n\geq 2$:
$$M_n=\psi\Big(
\frac{\Vert f_1\Vert_\infty+\cdots+\Vert f_{n-1}\Vert_\infty}{ 2}\Big)\,,$$
then $\big(M_n+\psi(x_n)\big)\beta_n\leq 1/2^n$;\par
\item[(iii)] for every $n\geq 1$, $\Vert g_n\Vert_\psi\geq 1/2$, with 
$g_n=f_n\,\ind_{\{\vert f_n\vert>\alpha_n\}}$.
\end{itemize}

For this, we start with $\beta_1$ such that $\big(1+\psi(x_1)\big)\beta_1=1/2$. Since the norms 
$\Vert\ \Vert_\psi$ and $\Vert\ \Vert_2$ are not equivalent on $X_0$, there is an $f_1\in X_0$ with 
$\Vert f_1\Vert_\psi=1$ and $\P (\{\vert f_1\vert>\alpha_1\})\leq \beta_1$. Suppose now that 
$f_1,\ldots,f_{n-1}$ and $\beta_1,\ldots,\beta_{n-1}$ have been constructed. We choose then $\beta_n\leq 1/2^n$ 
in order that $\big(M_n+\psi(x_n)\big)\beta_n\leq 1/2^n$. Since the norms $\Vert\ \Vert_\psi$ and 
$\Vert\ \Vert_2$ are not equivalent on $X_0$, we can find $f_n\in X_0$ such that $\Vert f_n\Vert_\psi=1$ and
$\Vert f_n\Vert_2$ is so small that
$$\P (\{\vert f_n\vert>\alpha_n\})\leq \beta_n\,.$$
Since $\Vert f_n-g_n\Vert_\psi\leq (1/a)\,\Vert f_n-g_n\Vert_\infty \leq \alpha_n/a$,
we have $\Vert g_n\Vert_\psi\geq \Vert f_n\Vert_\psi-\alpha_n/a\geq 1/2$,
and that finishes the construction.\par
Now, consider
$$g=\sum_{n=1}^{+\infty} \vert g_n\vert\,.$$
Set $A_n=\{\vert f_n\vert>\alpha_n\}$ and, for $n\geq 1$:
$$B_n=A_n\setminus \bigcup_{j>n} A_j\,.$$
We have $\P \big(\limsup A_n\big)=0$, because 
$\sum_{n\geq 1} \P (A_n)\leq \sum_{n\geq 1} \beta_n<+\infty$. Now $g$ vanishes out of 
$\bigcup_{n\geq 1} B_n\cup \big(\limsup A_n\big)$ and 
$\int_{B_n} \psi(\vert g_n\vert)\,d\P \leq \int_\Omega \psi(\vert f_n\vert)\,d\P \leq 1$. Therefore

\begin{align*}
\int_\Omega \psi\Big(\frac{\vert g\vert}{4}\Big)\,d\P
&= \sum_{n=1}^{+\infty} \int_{B_n} \psi\Big(\frac{\vert g\vert}{4}\Big)\,d\P \\
&\leq \sum_{n=1}^{+\infty} \int_{B_n} \frac{1}{ 2}\bigg[\psi\Big(
\frac{\Vert f_1\Vert_\infty+\cdots+\Vert f_{n-1}\Vert_\infty}{ 2}\Big)+
\psi\Big(\frac{\vert g_n\vert}{ 2}\Big)\bigg]\,d\P \\
&\hskip 5mm
\hbox{by convexity of $\psi$ and because $g_j=0$ on $B_n$ for $j>n$}\\
&\leq \frac{1}{ 2}\sum_{n=1}^{+\infty} M_n\,\P (A_n)
+ \frac{1}{ 2}\sum_{n=1}^{+\infty} \frac{1}{2^n} \int_{B_n} \! \! \psi(\vert g_n\vert)\,d\P +
\frac{1}{2}\sum_{n=1}^{+\infty} \psi(x_n)\,\P (A_n) \\
&\leq \frac{1}{2}\sum_{n=1}^{+\infty} \big(M_n+\psi(x_n)\big)\beta_n +
\frac{1}{2} \sum_{n=1}^{+\infty} \frac{1\,}{ 2^n} \leq 1\,.
\end{align*}
Hence $g\in L^{\psi}(\Omega)$.\par\vskip 1mm

It follows that the series ${\sum_{n\geq 1} g_n}_{\phantom j}$ is weakly unconditionally Cauchy in $X$. Since 
$\Vert g_n\Vert_\psi\geq 1/2$, it has, by the Bessaga-Pe\l czy\'nski's theorem, a subsequence which is equivalent
to the canonical basis of $c_0$. The same is true for $(f_n)_{n\geq 1}$ since
$$\sum_{n=1}^{+\infty} \Vert f_n - g_n\Vert_\psi
\leq \frac{1}{ a}\sum_{n=1}^{+\infty} \Vert f_n - g_n\Vert_\infty\leq
\frac{1}{ a}\sum_{n=1}^{+\infty} \alpha_n < 1\,.$$
That ends the proof.\hfill\qed
\medskip

Of course, the proof shows that the assumption that the norm $\Vert\ \Vert_\psi$ is not equivalent to 
$\Vert\ \Vert_2$ can be replaced by the non-equivalence of $\Vert\ \Vert_\psi$ with many other norms. We only
used the fact that the topology of convergence in measure is not equivalent on $X_0$ to the topology defined by 
$\Vert\ \Vert_\psi$.\par
\smallskip
When we apply this result to the probability space $(G,m)$, we get (see \cite{LLQR}, Theorem 2.3):
\goodbreak

\begin{theoreme}
Let $\psi$ be as in Theorem~\ref{theo 1} and let $G$ be a compact abelian group. Then, for 
$\Lambda\subseteq \Gamma=\hat G$, either $L^\psi_\Lambda$ has cotype 2, or it contains $c_0$.\par
In particular, either $\Lambda$ is a Sidon set and $L^{\Psi_2}_\Lambda=L^2_\Lambda$, or 
$L^{\Psi_2}_\Lambda$ contains $c_0$ (and so it has no finite cotype).
\end{theoreme}

\noindent{\bf Proof.} Observe that when $L^\psi_\Lambda \ne L^2_\Lambda$, the norms $\Vert\ \Vert_\psi$ 
and $\Vert\ \Vert_2$ are not equivalent on $X_0={\cal P}_\Lambda$, the subspace of the trigonometric 
polynomials whose spectrum is contained in $\Lambda$. So the first part follows directly from
Theorem~\ref{theo 1}. The second one follows from Pisier's characterization of Sidon sets (\cite{Pisier}, 
Th\'eor\`eme 6.2): $\Lambda$ is a Sidon set if and only if $L^{\Psi_2}_\Lambda=L^2_\Lambda$.\hfill\qed
\medskip

\noindent{\bf Remark.} It is proved in \cite{LLQR}, Theorem 2.3, that $\Lambda$ is a $\Lambda(\psi)$-set 
({\it i.e.} $L^\psi_\Lambda= L^2_\Lambda$) when 
$L^{\Psi_2}_\Lambda\subseteq L^\psi_\Lambda\subseteq L^2_\Lambda$ and $L^\psi_\Lambda$ has cotype 2.
\bigskip\goodbreak

\section{Uniform convergence}

A function $f\in {\cal C}(\T)$ is said to have a {\it uniformly convergent Fourier series} if 
$\Vert S_k(f) - f\Vert_\infty\mathop{\longrightarrow}\limits_{k\to+\infty} 0$, where
$$S_k(f)=\sum_{j=-k}^k \hat f(j)\,e_j\,.$$
\par
The space $U(\T)$ of {\it uniformly convergent Fourier series} is the space of all such $f\in {\cal C}(\T)$. With the 
norm
$$\Vert f\Vert_U=\sup_{k\geq 1} \Vert S_k(f)\Vert_\infty\,,$$
$U(\T)$ becomes a Banach space.\par

A set $\Lambda\subseteq \Z$ is said to be a {\it set of uniform convergence} ({\it $UC$-set}) if 
$U_\Lambda={\cal C}_\Lambda$ as linear spaces. They are then isomorphic as Banach spaces. There exist sets 
$\Lambda$ which are not $UC$-sets, but for which ${\cal C}_\Lambda$ does not contain $c_0$ (for instance, a 
Rosenthal set which contains arbitrarily long arithmetical progressions \cite{Rosenthal}). For $U_\Lambda$ the 
situation is different; we have:

\begin{theoreme}\label{theo 3.1}
If $\Lambda$ is not a $UC$-set, then $U_\Lambda$ contains $c_0$.
\end{theoreme}

\begin{corollaire}
If $U_\Lambda$ has a finite cotype, then $\Lambda$ is a Sidon set.
\end{corollaire}

\noindent{\bf Proof.} If $U_\Lambda$ has a finite cotype, it cannot contain $c_0$. Hence $U_\Lambda$ is 
isomorphic to ${\cal C}_\Lambda$. It follows that ${\cal C}_\Lambda$ has a finite cotype, and so $\Lambda$ is a 
Sidon set, by Bourgain-Milman's theorem \cite{Bourgain-Milman}.\hfill\qed
\medskip

\noindent{\bf Remark.} This result was proved in \cite{LLQR}, Theorem 1.2, by adapting the proof of Bourgain and 
Milman. Now it becomes clear why this proof happened to mimic the original one.
\medskip

\noindent{\bf Proof of Theorem~\ref{theo 3.1}.} Since $\Lambda$ is not a $UC$-set, there exists a trigonometric 
polynomial $P_1\in {\cal C}_\Lambda$ such that $\Vert P_1\Vert_U=1$ and $\Vert P_1\Vert_\infty\leq 1/2$. 
Let $N_1\geq 2$ such that $\hat{P_1}(n)=0$ for $\vert n\vert\geq N_1$. The spaces 
$U_{\Lambda\setminus \Lambda\cap \{-N_1+1,\ldots, 0,\ldots, N_1-1\}}$ and
${\cal C}_{\Lambda\setminus \Lambda\cap \{-N_1+1,\ldots, 0,\ldots, N_1-1\}}$ 
remain non-isomorphic, and so we can find a trigonometric polynomial $P_2$ such that $\hat{P_2}(n)=0$ for 
$\vert n\vert\leq N_1-1$ with $\Vert P_2\Vert_U=1$ and $\Vert P_2\Vert_\infty\leq 1/4$. Carrying on this 
construction, we get a sequence of integers $2\leq N_1< N_2<\cdots$ and a sequence of trigonometric polynomials 
$P_l\in {\cal C}_\Lambda$ such that $\Vert P_l\Vert_U=1$, $\Vert P_l\Vert_\infty\leq 1/2^l$ and
$\hat{P_l}(n)=0$ for $n\notin \{\pm N_{l-1}, \ldots, \pm (N_l -1)\}$.\par
Now, fix an integer $L\geq 1$ and a sequence $a_1,\ldots, a_L$ of complex numbers. For each $k\geq 1$, let $l_k$ 
such that $N_{l_k}\leq k <N_{l_k +1}$. We have, when $L\geq l_k +1$:
\begin{align*}
\Big\Vert S_k\Big(\sum_{l=1}^L a_l P_l\Big)\Big\Vert_\infty
&\leq \Big\Vert \sum_{l=1}^{l_k} a_l P_l\Big\Vert_\infty
+ \Vert a_{l_k+1} S_k (P_{l_k+1})\Vert_\infty \\
&\leq \max_{1\leq j\leq l_k} \vert a_j\vert\hskip 1mm
\sum_{l=1}^{l_k} \Vert P_l\Vert_\infty + \vert a_{l_k+1}\vert\,\Vert P_{l_k+1}\Vert_U\\
&\leq 2\,\max \{\vert a_1\vert,\ldots, \vert a_{l_k}\vert, \vert a_{l_k+1}\vert, \ldots, \vert a_L\vert\}\,.
\end{align*}
The inequality
$\big\Vert S_k\big(\sum_{l=1}^L a_l P_l\big)\big\Vert_\infty \leq
2\,\max \{\vert a_1\vert,\ldots, \vert a_{l_k}\vert, \vert a_{l_k+1}\vert, \ldots, \vert a_L\vert\}$
remains trivially true for $L\leq l_k$, because in this case
$S_k\big(\sum_{l=1}^L a_l P_l\big)=\sum_{l=1}^L a_l P_l$. Therefore we get
$$\Big\Vert \sum_{l=1}^L a_l P_l\Big\Vert_U\leq
2\,\max\{\vert a_1\vert,\ldots, \vert a_L\vert\}\,.$$
It follows that the series $\sum_{l\geq 1} P_l$ is weakly unconditionally Cauchy. Since it is obviously not convergent, 
$U_\Lambda$ contains a subspace isomorphic to $c_0$ by Bessaga-Pe\l czy\'nski's theorem (see 
\cite{Diestel book}, pages 44--45, Theorem 6 and Theorem 8).\hfill\qed
\medskip

\noindent{\bf Remark 1.} There is the stronger notion of $CUC$-set. $\Lambda\subseteq \Z$ is a $CUC$-set if 
$\Big\Vert\sum_{j=k_1}^{k_2} \hat f(j)\,e_j -f\Big\Vert_\infty
\mathop{\longrightarrow}\limits_{\mathop{\scriptstyle k_1\to -\infty}\limits_{\scriptstyle k_2\to +\infty}} 0$ 
for every $f\in {\cal C}_\Lambda$. Obviously, for subsets of $\N$, the two notions coincide. Theorem~\ref{theo 3.1} 
is not valid for $CUC$-sets: let $H$ be a Hadamard lacunary sequence. Then $\Lambda=H-H$ is not a $CUC$-set 
(Fournier \cite{Fournier}), but it is $UC$ and Rosenthal, so that $U_\Lambda={\cal C}_\Lambda$ does not
contain $c_0$.\par
However, it is not known whether ${\cal C}_{\Lambda_1\cup \Lambda_2}$ lacks $c_0$ whenever this is true for 
${\cal C}_{\Lambda_1}$ and ${\cal C}_{\Lambda_2}$. If we replace the space ${\cal C}(G)$ by $U(\T)$, the 
answer is in the negative. Indeed, J. Fournier shows (\cite{Fournier}), completing Soardi and Travaglini's work 
\cite{S-T}, that there exist two $UC$-sets $\Lambda_1, \Lambda_2\subseteq \Z$, which are Rosenthal sets, but
$\Lambda_1\cup \Lambda_2= H+H-H$ is not $UC$. Therefore $U_{\Lambda_1}={\cal C}_{\Lambda_1}$ and
$U_{\Lambda_2}={\cal C}_{\Lambda_2}$ do not contain $c_0$, though $U_{\Lambda_1\cup \Lambda_2}$ 
contains $c_0$.
\medskip

\noindent{\bf Remark 2.} $UC$-sets $\Lambda$ for which ${\cal C}_\Lambda$ contains $c_0$ are constructed 
in \cite{LQR}.
\medskip

\noindent{\bf Remark 3.} We stated Theorem~\ref{theo 3.1} for uniform convergence because it is the classical case. 
Actually, J. Fournier (\cite{Fournier}, page 72) and S. Hartman (\cite{Hartman}, page 107) introduced the space 
$L^1-UC$ as the set of all $f\in L^1(\T)$ for which 
$\Vert S_k(f) -f\Vert_1\mathop{\longrightarrow}\limits_{k\to+\infty} 0$. It is normed by 
$\Vert f\Vert_{UL^1}=\sup_{k\geq 1} \Vert S_k(f)\Vert_1$. $\Lambda$ is said to be an $L^1-UC$-set if 
$(L^1-UC)_\Lambda=L^1_\Lambda$. The same proof as above shows that if 
$(L^1-UC)_\Lambda \not= L^1_\Lambda$, then $(L^1-UC)_\Lambda$ contains $c_0$. More generally, let 
$\Lambda\subseteq \Z$ and let $X$ be a Banach space contained, as a linear subspace, in $L^1(\T)$ such that the
linear space generated by $X\cap \Lambda$ is dense in $X$. We can define $X-UC$ in an obvious way, and we 
have: if $X-UC$ is not isomorphic to $X$, then it contains $c_0$.
\medskip

We give another consequence of Theorem~\ref{theo 3.1}. Recall (see \cite{Meyer}) that 
$\Lambda\subseteq \Gamma$ is a \emph{Riesz set} if every measure with spectrum in $\Lambda$ is absolutely 
continuous, with respect to the Haar measure (in short, ${\cal M}_\Lambda=L^1_\Lambda$).

\begin{corollaire}
If $U_\Lambda$ does not contain $c_0$, then $\Lambda$ is a Riesz set.
\end{corollaire}

\noindent{\bf Proof.} If $U_\Lambda\not\supseteq c_0$, then $U_\Lambda={\cal C}_\Lambda$, by 
Theorem~\ref{theo 3.1}, and so ${\cal C}_\Lambda\not\supseteq c_0$. It follows then that $\Lambda$ is a
Riesz set (F. Lust-Piquard \cite{LP1}, her first Th\'eor\`eme 3.1). Let us recall why. For 
$\mu\in {\cal M}_\Lambda$, the convolution operator 
$C_\mu\colon f\in {\cal C}(G)\mapsto f\ast\mu\in {\cal C}_\Lambda\subseteq {\cal C}(G)$
is weakly compact, because ${\cal C}(G)$ has Pe{\l}czy\'nski's property $(V)$ and 
${\cal C}_\Lambda\not\supseteq c_0$. Its adjoint operator 
$\nu\in {\cal M}(G)\mapsto \nu\ast\mu\in {\cal M}_\Lambda$ is also weakly compact. Hence, if $(K_j)_j$ is an 
approximate unit for the convolution, there is a sequence $(j_n)_n$ such that $K_{j_n}\ast\mu$ is weakly
convergent. Since $K_j\ast \mu$ converges weak-star to $\mu$, it follows that $\mu\in L^1_\Lambda$.\hfill\qed
\medskip

\noindent{\bf Remark.} Another proof can be given, without using Theorem~\ref{theo 3.1}, but using that 
$U(\T)$ has Pe{\l}czy\'nski's property $(V)$ (Saccone \cite{Saccone2}, Theorem 2.2; for $U_\N(\T)$, see 
Bourgain \cite{Bourgain1}, Lemme 2 and Lemme 3, and Saccone \cite{Saccone1}, Theorem 4.1). Then, as before,
$K_{j_n}\ast\mu$ is weakly convergent, in $U(\T)^\ast$ this time. So there are convex combinations which 
converge in the norm of $U(\T)^\ast$. But then they converge in the norm of $U_\N(\T)^\ast$, and so 
$u\in L^1(G)$ (see D. Oberlin \cite{Oberlin}, page 310). Note that Oberlin's argument (as well as Bourgain's one) 
depends on Carleson's Theorem ({\it via} \cite{Vi}).
\bigskip

\section{Invariant means and Hilbert sets}

An \emph{invariant mean} $M$ on $L^\infty(G)$ is a continuous linear functional on $L^\infty(G)$ such that 
$M(\ind)=\Vert M\Vert=1$ and $M(f_x)=M(f)$ for every $f\in L^\infty(G)$. The Haar measure $m$ defines an 
invariant mean, and W. Rudin (\cite{Rudin}) showed that, for infinite compact abelian groups $G$, there always 
exist other invariant means on $L^\infty(G)$. A function $f\in L^\infty(G)$ has a unique invariant mean if 
$M(f)=\hat f(0)$ for every invariant mean $M$ on $L^\infty(G)$. Every continuous function (or, even, every 
Riemann-integrable function: \cite{Rubel-Shields}, page 38, or \cite{Ta1}) has a unique invariant mean.

\begin{definition}
A subset $\Lambda$ of $\Gamma=\hat G$ is called a \emph{Lust-Piquard set} if $\gamma f$ has a unique 
invariant mean for every $f\in L^\infty_\Lambda$ and every $\gamma\in \Gamma$.
\end{definition}

In other words, $\Lambda$ is a Lust-Piquard set if for every invariant mean $M$ on $L^\infty(G)$ and every 
$f\in L^\infty_\Lambda$, one has: 
$$M(\gamma f)=\hat f(-\gamma).$$

In \cite{LP2} (and then in \cite{Li1}; see also \cite{LP4}), F. Lust-Piquard called them {\it totally ergodic sets}. 
We use a different name because J. Bourgain (\cite{Bourgain2}, 2.I, page 206), used the terminology ``ergodic set'' 
for another property (see also \cite{LQR}).\par
Note that it is required that the invariant means agree on 
$\bigcup_{\gamma\in \Gamma} L^\infty_{\Lambda-\gamma}$, and not only on $L^\infty_\Lambda$, because 
the invariant means may coincide on $L^\infty_\Lambda$ for trivial reasons; for instance, all the invariant 
means are equal to $0$ on $L^\infty_{2\Z+1}$ (since $f(x+1/2)=-f(x)$ for $f\in L^\infty_{2\Z+1}$). It is clear that 
if $\Lambda$ is a Lust-Piquard set, then $\Lambda-\gamma$ is also a Lust-Piquard set for every 
$\gamma\in\Gamma$.\par
It is obvious that every Rosenthal set is a Lust-Piquard set (since every continuous function has a unique invariant 
mean), and it is shown in \cite{Li1} that every Lust-Piquard set is a Riesz set. On the other hand, Y. Katznelson
(see \cite{Rubel-Shields}, pages 37--38) proved that $\N$ is not a Lust-Piquard set.
\medskip

F. Lust-Piquard (\cite{LP3}, Theorem 2 and Theorem 4) showed that $\Lambda={\P}\cap (5\Z+2)$, where ${\P}$ 
is the set of the prime numbers, is totally ergodic (a Lust-Piquard set, with our terminology) although 
${\cal C}_\Lambda$ contains $c_0$.\par
In the following theorem, we give another example of such a situation. Let us recall that $H\subseteq \Z$ is a 
\emph{Hilbert set} if there exist two sequences of integers $(p_n)_{n\geq 1}$ and $(q_n)_{n\geq 1}$, with
$q_n\not=0$, such that
$$H=\bigcup_{n\geq 1} \big\{p_n+\sum_{k=1}^n \eps_k\,q_k\,;\ \eps_1,\ldots,\eps_n= 0\ {\rm or}\ 1\big \}.$$
It is shown in \cite{Li2}, Theorem 2, that ${\cal C}_H$ contains $c_0$ when $H$ is a Hilbert set.

\begin{theoreme}\label{theo 4.2}
There exists a Hilbert set $H\subseteq \N$ which is a Lust-Piquard set.
\end{theoreme}

We begin with a lemma, which is implicit in \cite{LP3}, proof of Theorem 4.

\begin{lemme}\label{lemme 4.3}
The family of Lust-Piquard sets in $\Gamma$ is localizable for the Bohr topology.
\end{lemme}

Let us recall that the \emph{Bohr topology} of a discrete abelian group $\Gamma$ is the topology of pointwise 
convergence, when $\Gamma$ is seen as a subset of ${\cal C}(G)$; it is also the natural topology on $\Gamma$ 
as a subset of the dual group of $G_d$, the group $G$ with the discrete topology. A class ${\cal F}$ of subsets of 
$\Gamma$ is \emph{localizable for the Bohr topology} if $\Lambda\in {\cal F}$ whenever for every 
$\gamma\in \Gamma$ there is a neighbourhood $V_\gamma$ of $\gamma$ for the Bohr topology such that
$\Lambda\cap V_\gamma\in {\cal F}$. This notion is due to Y. Meyer (\cite{Meyer}).
\par
For the sake of completeness, we shall give a proof.
\medskip

\noindent{\bf Proof of Lemma~\ref{lemme 4.3}.} We are going to prove that if $V_\gamma$ is a neighbourhood 
of $\gamma\in \Gamma$ such that $\Lambda\cap V_\gamma$ is a Lust-Piquard set, then $\overline\gamma f$ 
has a unique invariant mean for every $f\in L^\infty_\Lambda$, and that will prove the lemma.\par
By the regularity of the algebra $L^1(G_d)=\ell_1(G)={\cal M}_d(G)$, there exists a discrete measure 
$\nu\in {\cal M}_d(G)$ such that $\hat \nu(\gamma)=1$ and $\hat\nu=0$ outside $V_\gamma$. Since
$(\overline \gamma f)\ast (\overline \gamma \nu)\in L^\infty_{(\Lambda\cap V_\gamma)-\gamma}$, and since
$(\Lambda\cap V_\gamma)-\gamma$ is a Lust-Piquard set, we have:
$$M\big((\overline \gamma f)\ast (\overline \gamma \nu)\big)=
\widehat{[(\overline \gamma f)\ast (\overline \gamma \nu)]^{}} (0)=
\hat f(\gamma)\,\hat\nu(\gamma)=\hat f(\gamma).$$
But $\overline\gamma\nu$ is a discrete measure, and we have, for every discrete measure $\mu$:
$$M(\mu\ast g)=M(g)\,\hat\mu(0)$$
for every $g\in L^\infty(G)$ and every invariant mean $M$. This is so since, if $\mu=\sum_k a_k\,\delta_{x_k}$, 
with $\sum_k\vert a_k\vert<+\infty$, we have 
$$M(\mu\ast g)=\sum_k a_k\,M(g_{x_k})= \sum_k a_k\,M(g).$$
\par
Hence $M(\overline \gamma f)=\hat f(\gamma)$, as required.\hfill\qed
\medskip

\noindent{\bf Proof of Theorem~\ref{theo 4.2}.} We are going to construct a Hilbert set $H\subseteq \N$ which is 
discrete in $\Z$ for the Bohr topology. For such a set, there is, for every $k\in\Z$, some Bohr-neighbourhood 
$V_k$ of $k$ such that $H\cap V_k$ is finite. Therefore, we have $L^\infty_{H\cap V_k}={\cal C}_{H\cap V_k}$, 
and so $H\cap V_k$ is a Lust-Piquard set.\par
Let $(d_n)_{n\geq 0}$ be a strictly increasing sequence of positive integers such that:
$$d_n\mid d_{n+1}\,,\quad n\geq 0\,, \qquad  \sum_{n=0}^{+\infty} \frac{2^{n+1}}{d_n} <1\,.$$
\par
For every $k\in\Z$, consider:
$$V(k)=k+d_{\vert k\vert}\Z\,,$$
which is a Bohr-neighbourhood of $k$.\par
Now, we are going to show that we can choose, for every $n\geq 0$, an
integer $r_n\in \{0,1,2,\ldots,d_n-1\}$ such that, if
$$H_n= d_n + r_n +\Big\{\sum_{l=0}^{n-1} \eps_l d_l \,;\ \eps_l=0\ {\rm or}\ 1\Big\}\,,$$
then
$$H_p\cap V(k)=\emptyset$$
for every $k\in \Z$ and every $p>\vert k\vert$. The set $H=\bigcup_{n\geq 0} H_n$ will be the required set.\par
\vskip 1mm

We are going to do this by induction. First, we may choose an arbitrary $r_0\in \{0,1,2,\ldots,d_0-1\}$, and we set 
$H_0=\{d_0+r_0\}$. Suppose now that we have found $r_1, r_2,\ldots, r_{p-1}$ such that the previous
conditions are fulfilled:
$$\hskip 4cm H_j\cap V(k)=\emptyset\,,\hskip 5mm{\rm for} \ 1\leq j\leq p-1\,,\ \vert k\vert<j\,.$$
To find $r_p$, note that $m\in H_p\cap V(k)$ if and only if
\begin{equation}
m\in k+ d_{\vert k\vert}\Z \tag*{(1)}
\end{equation}
and there exist $\eps_0,\eps_1,\ldots,\eps_{p-1}\in \{0,1\}$ such that
\begin{equation}
m=d_p+r_p+\sum_{l=0}^{p-1} \eps_l d_l\,.\tag*{(2)}
\end{equation}
Since, for $0\leq l<p$, one has $d_l\mid d_{l+1}\mid\cdots\mid d_p$, conditions $(1)$ and $(2)$ are equivalent 
to $r_p\equiv 0\ ({\rm mod}\,d_0)$, for $k=0$, and, for $1\leq l=\vert k\vert<p$, to:
$$\hskip 1cm k \equiv r_p+\sum_{j=0}^{\vert k\vert-1} \eps_j d_j\hskip 5mm ({\rm mod}\,d_{\vert k\vert})\,.$$
For each such $k$ ($0\leq \vert k\vert<p$), there are
$$\frac{d_p}{d_{\vert k\vert}}\,\cdot 2^{\vert k\vert}$$
possible choices for $r_p$. As
$$\frac{d_p}{d_0}+ 2\sum_{l=1}^{p-1} 2^l \frac{d_p}{d_l}
\leq \frac{d_p}{d_0}+2\sum_{l=1}^{+\infty} 2^l \frac{d_p}{d_l} < d_p\,,$$
by hypothesis, we can find an $r_p\in\{0,1,\ldots, d_p-1\}$ such that the set $H_p$ constructed from it verifies
$H_p\cap V(k)=\emptyset$ for $\vert k\vert<p$. That ends the proof.
\hfill\qed
\medskip

\noindent{\bf Remark 1.} Some particular Hilbert sets are the $IP$-sets, {\it i.e.} the sets $F$ for which there exists 
a sequence $(p_n)_{n\geq 1}$ of integers such that
$$F=\Big\{\sum_{k=1}^n \eps_k\,p_k\,;\ \eps_1,\ldots,\eps_n=0\ {\rm or}\ 1 \,,\ n\geq 1\Big\}\,.$$\par

\noindent{\bf Question.} Does there exist an $IP$-set $F$ which is a Lust-Piquard set?\par\bigskip

Every point of an $IP$-set $F$ is non-isolated in $F$ (see \cite{Furst}, Theorem 2.19; note that every point of an 
$IP$-set is inside the translation by this point of a sub-$IP$-set). Therefore we cannot use an argument similar to 
that of the previous theorem. Hilbert sets and $IP$-sets are different in several ways. For instance, every set 
$\Lambda\subseteq \Z$ which has a positive uniform density contains a Hilbert set (\cite{Hindman},
Theorem 11.11; \cite{Li2}, Theorem 4), but not necessarily an $IP$-set (\cite{Hindman}, Theorem 11.6; 
\cite{Nath}, page 151). Another difference is that ${\cal C}_\Lambda$ never has the Unconditional Metric 
Approximation Property if $\Lambda\subseteq \Z$ is an $IP$-set ([23], Proposition 11), but can have
this property when $\Lambda$ is a Hilbert set (\cite{Li3}, Theorem 10).\par
\smallskip

\noindent{\bf Remark 2.} Let ${\cal F}$ be a class of subsets of $\Gamma$, which contains all the finite sets, and 
which is localizable for the Bohr topology. It follows from the proof of Theorem~\ref{theo 4.2} that such a class must
contain some Hilbert sets. In particular ${\cal F}$ has to contain sets $\Lambda$ such that $\Lambda$ contains 
parallelepipeds of arbitrarily large dimensions. Note that this last assertion is actually implicit in \cite{LP3}.
Indeed, by Dirichlet's theorem, $\sum_{n\in {\P}\cap(5\Z+2)} \frac{1}{n}=+\infty$, and
by \cite{Miheev}, Corollary 2, we have $\sum_{n\in \Lambda} \frac{1}{n}<+\infty$ when $\Lambda$ does not 
contain parallelepipeds of arbitrarily large dimensions. It is known that the sets belonging to the following classes
cannot contain parallelepipeds of arbitrarily large dimensions:\par\smallskip

\begin{itemize}
\item [a)] $\Lambda(p)$-sets (see \cite{Miheev}, Theorem 3, and \cite{Fournier-Pigno}, Theorem 4);
\item [b)] $UC$-sets (\cite{Fournier-Pigno}, Theorem 4);
\item [c)] $p$-Sidon sets (\cite{J-W}, Lemma 1);\par
\item [d)] stationary sets (\cite{Lef}, Proposition 2.5);\par
\item [e)] $q$-Rider sets (see \cite{LQR} or \cite{LLQR} for the definition). Note that, for $1\leq q<4/3$, 
$q$-Rider sets are  $p$-Sidon sets, for every $p>q/(2-q)$ (see \cite{Lef-Ro}), and so the result is in c). For 
$4/3\leq q<2$, there is no explicit published proof of that, and therefore we shall give one in 
Proposition~\ref{prop 4.4}, after this Remark.
\end{itemize}

Hence these classes are not localizable for the Bohr topology.\par
This remark shows that there is no hope to construct sets of the above classes by way of localization.

\begin{proposition}\label{prop 4.4}
If $\Lambda$ is a $q$-Rider set, $1\leq q<2$, then $\Lambda$ cannot contain parallelepipeds of arbitrarily large 
dimensions.
\end{proposition}

\noindent{\bf Proof.} A Sidon set (with constant less than $10$, say) inside a parallelepiped $P$ of size $2^n$ 
cannot contain more than $C\,n\log n$ elements (\cite{Kahane}, Chapter 6, \S\kern 1mm 3, Theorem 5, page 71),
whereas if $P$ were contained in a $q$-Rider set, it should contain a quasi-independent (hence Sidon with constant 
less than $10$) set of size at least $C\,2^{\eps n}$, with $\eps=(2-q)/q$ (\cite{Luis}, or \cite{Luis-T}, Teorema 2.3).
\hfill\qed
\medskip

Note that another proof of Proposition~\ref{prop 4.4} is implicit in \cite{J-W}. Indeed the proof given in 
\cite{J-W}, Lemma\kern 1mm 1, that $p$-Sidon sets share this property only uses the fact, proved in 
\cite{Bo-Pyt}, Eq. $(9)$, that if $\Lambda$ is a $p$-Sidon set, then, with $\alpha=2p/(3p-2)$, there is a constant 
$C>0$ such that $\Vert f\Vert_r\leq C\,\sqrt r\,\Vert \hat f\Vert_\alpha$ for all $r\geq 2$ (equivalently:
$\Vert f\Vert_{\Psi_2}\leq C'\,\Vert \hat f\Vert_\alpha$) for every $f\in {\cal C}_\Lambda$. Now the 
fourth-named author proved that these inequalities characterize $p$-Rider sets (\cite{Luis}; see also 
\cite{Luis-T}, Teorema 2.3).\par

\goodbreak

\section{Complemented subspaces}

Since $\Lambda$ is a Rosenthal set if $L^\infty_\Lambda={\cal C}_\Lambda$, it is natural to ask whether 
$\Lambda$ is a Rosenthal set if there exists a projection from $L^\infty_\Lambda$ onto ${\cal C}_\Lambda$. 
We have not been able to answer this, even if this projection were to have norm $1$ (see \cite{Godefroy}, where the 
condition that the space does not contain $\ell_1$ is crucial), but we have a partial result. Recall that it is not known 
whether ${\cal C}_\Lambda\not \supseteq c_0$ implies that $\Lambda$ is a Rosenthal
set.\par

\begin{theoreme}
Let \hfill $\Lambda\subseteq \Gamma$ \hfill be \hfill such \hfill that \hfill there \hfill exists \hfill a 
\hfill surjective \hfill projection\\
$P\colon L^\infty_\Lambda\to {\cal C}_\Lambda$. Then ${\cal C}_\Lambda$ does not contain $c_0$. 
Moreover, every Riemann-integrable function in $L^\infty_\Lambda$ is actually in ${\cal C}_\Lambda$.
\end{theoreme}

Recall that a function $h\colon G\to \C$ is \emph{Riemann-integrable} if it is bounded and almost everywhere 
continuous. Actually, the last assertion of the proposition means that every element of $L^\infty_\Lambda$ which 
contains a Riemann-integrable function contains also a continuous one.
\smallskip

\noindent{\bf Proof.} 1) By \cite{Li2}, Proposition 14, if ${\cal C}_\Lambda$ contains $c_0$, there is a sequence 
$(f_n)_{n\geq 1}$ in ${\cal C}_\Lambda$,  which is equivalent to the canonical basis of $c_0$, and whose 
$w^\ast$ linear span $F$ in $L^\infty_\Lambda$ is isomorphic to $\ell_\infty$. The restriction $P_{\mid F}$ 
is a projection from $F$ onto a subspace of ${\cal C}_\Lambda$ which contains 
$E=\overline{\phantom{i}\hskip -0,3mm{\rm span}}\,\{f_n\,;\ n\geq 1\}$.\par
Observe that $E$ is a separable subspace of ${\cal C}_\Lambda$. So there exists a countable subset 
$\Lambda_1\subseteq \Lambda$ such that $E\subseteq {\cal C}_{\Lambda_1}$. Moreover, there exists a 
countable subgroup $\Gamma_0\subseteq \Gamma$ such that $\Lambda_1$ is contained in $\Gamma_0$. 
Taking $\Lambda_0=\Lambda\cap \Gamma_0$, we have $E\subseteq {\cal C}_{\Lambda_0}$, and 
${\cal C}_{\Lambda_0}$ is a separable space.\par
The set $\Gamma_0$ being a subgroup, there exists a measure $\mu$ on $G$ whose Fourier transform is 
$\hat\mu = \ind_{\Gamma_0}$. The map $f\mapsto f\ast \mu$ gives a projection from ${\cal C}_{\Lambda}$ 
onto ${\cal C}_{\Lambda_0}$, and Sobczyk's theorem gives a projection from ${\cal C}_{\Lambda_0}$ onto $E$. 
So there exists a projection from $F\simeq \ell_\infty$ onto $E\simeq c_0$, which is a contradiction.\par
\smallskip
2) We first assume that the group $G$ is metrizable, so that ${\cal C}(G)$ is separable. \par
Let $RI_\Lambda$ be the subspace of $L^\infty_\Lambda$ consisting of Riemann-integrable functions (more 
precisely: the elements of $L^\infty_\Lambda$ which have a Riemann-integrable representative).\par
Consider the restriction of $P$ to $RI_\Lambda$. For $f\in RI$, the set
$$\{x\mapsto \xi(f_x)\,;\ \xi\in L^\infty(G)^\ast,\ \Vert \xi\Vert\leq 1\}$$
is stable (\cite{Ta3}, Theorem 15-6\kern 0,5mm c)). Let $\mu\in ({\cal C}_\Lambda)^\ast$, and set
$\phi(x,y)=\big(P^\ast\mu_y)(f_x)$ for $x,y\in G$. The map $x\in G\mapsto f_x\in L^\infty(G)$ is scalarly 
measurable (\cite{Ta2}, Theorem 16) and $y\mapsto P^\ast\mu_y$ is continuous for the $w^\ast$-topology. 
Moreover $\{x\mapsto (P^\ast\mu_y)(f_x)\,;\ y\in G\}$ is stable, so by \cite{Ta3}, Theorem 10-2-1, $\phi$ is 
measurable. Measurability refers here to the completion of the product measure $m\otimes m$ on $G\times G$, so 
in order to deduce that the map $x\in G\mapsto \phi(x,x)=(P^\ast \mu_x)(f_x)$ is measurable, we need the 
following lemma (note that our $\phi$ is bounded).

\begin{lemme}
Let $G$ be a metrizable compact abelian group, and $\phi\colon  G\times G\to \C$ a function such that:\par
1) $\phi\in  {\cal L}^\infty(G\times G)$;\par
2) the map $y\mapsto  \phi(x,y)$ is continuous, for every $x\in G$.\par
\noindent
Then the map $x\mapsto \phi(x,x)$ is measurable.
\end{lemme}

\noindent{\bf Proof.} $G$ being metrizable, there exists a bounded sequence $(f_n)_n$ in $L^1(G)$ such that
\begin{equation}
g(0)= \lim_{n\to \infty} \int_G f_n g \, dm \,, \qquad\hbox{for every $g\in {\cal C}(G)$.} \tag*{(3)}
\end{equation}
This sequence $(f_n)_n$ represents an approximate identity.\par
For every $n$, the function $(x,y)\mapsto f_n(x-y)\phi(x,y)$ is integrable in $G\times G$. Define
$$F_n(x)=\int_G f_n(x-y)\,\phi(x,y)\, dm(y) = \int_G f_n(t)\,\phi(x,x-t)\, dm(t).$$
By Fubini's theorem $F_n$ is defined almost everywhere, and is integrable. So $F_n$ is measurable, for every $n$. 
The lemma follows since, by $(3)$,
\begin{equation}
\phi(x,x)=\lim_{n\to\infty} F_n(x)\,, \qquad\hbox{ for every $x\in G$.}\tag*{$\square$}
\end{equation}

The fact that the map $x\in G\mapsto (P^\ast \mu_x)(f_x)=\,<\mu,\big[P(f_x)\big]_{-x}>$ is measurable means, 
since $\mu$ is arbitrary, that $x\mapsto \big[P(f_x)\big]_{-x}\in {\cal C}_\Lambda$ is scalarly measurable. Since 
we have assumed that ${\cal C}(G)$ is separable, this map is strongly measurable, by Pettis's measurability theorem
(\cite{Diestel-Uhl}, II \S\kern 0,5mm 1, Theorem 2). Now we showed in the beginning of the proof that 
${\cal C}_\Lambda$ does not contain $c_0$; so a result of J. Diestel \cite{Diestel} (see 
\cite{Diestel-Uhl}, II, \S\kern 0,5mm 3, Theorem 7) says that this map is Pettis-integrable, which means that if we 
define $Qf$ using 
$$<Qf,\mu>\,=\int_G <f_x,P^\ast(\mu_x)>\,dx\,,$$
for every $\mu\in ({\cal C}_\Lambda)^\ast$, then $Q$ maps $RI_\Lambda$ into ${\cal C}_\Lambda$, and not 
only into its bidual (see the definition of Pettis-integrability in \cite{Diestel-Uhl}, II, \S\kern 0,5mm 3, page 53, 
Definition 2, or in \cite{Ta3}, Definition 4-2-1).\par
Thus $Q$ is a projection from $RI_\Lambda$ onto ${\cal C}_\Lambda$ such that $Q(f_x)=(Qf)_x$ for every 
$f\in RI_\Lambda$ and every $x\in G$.\par
We want to prove that $Qf=f$ for every $f\in RI_\Lambda$, and for that we have to see that 
$\hat{Qf}(\gamma)=\hat f(\gamma)$ for every $\gamma\in \Gamma$. But it suffices to show that 
$\hat{Qf}(0)=\hat f(0)$, since, after replacing $\Lambda$ by $(\Lambda-\gamma)$ and $Q$ by
$Q_\gamma\colon L^\infty_{\Lambda-\gamma}\to {\cal C}_{\Lambda-\gamma}$, with 
$Q_\gamma(g)=\overline\gamma Q(\gamma g)$, we then get for $f\in RI_\Lambda$ with 
$g=\overline\gamma f$:
$$\hat{Qf}(\gamma)={[\overline\gamma (Qf)]}\hat {\hskip 2mm} (0) =
\hat{Q_\gamma g}(0)=\hat g(0)= \hat{(\overline\gamma f)}(0)=\hat f(\gamma)\,.$$
\par
So, let $f\in RI_\Lambda$. Every Riemann-integrable function has a unique invariant mean 
(\cite{Rubel-Shields}, Lemma 7; [44]); hence there are (\cite{Rubel-Shields}, Proposition page 38; or 
\cite{LP2}, Proposition 1) convex combinations $\sum\limits_{k\in I_n} c_{n,k}\, f_{x_{n,k}}$, $c_{n,k}>0$,
$\sum\limits_{k\in I_n} c_{n,k} =1$, of translates of $f$ which converge in norm to the constant function 
$\hat f(0)\,\ind$. We have:
$$Q\Big(\sum_{k\in I_n} c_{n,k} f_{x_{n,k}}\Big)\mathop{\longrightarrow}_{n\to+\infty}
Q\big[\hat f(0)\ind\big]=\hat f(0)\ind\,.$$
But $Q\big(\sum_{k\in I_n} c_{n,k} f_{x_{n,k}}\big)= \sum_{k\in I_n} c_{n,k} (Qf)_{x_{n,k}}$,
and its Fourier coefficient at $0$ is:
$$\sum_{k\in I_n} c_{n,k}\hat{Qf}(0) =\hat{Qf}(0)\,.$$
Therefore $\hat{Qf}(0)=\hat f(0)$.\hfill\qed
\par
\smallskip

3) In order to finish the proof, we have to explain why we may assume that $G$ is metrizable.\par
Let $\Lambda$ be as in the theorem, and $f\in RI_\Lambda$. As explained in the proof of the first part of the 
theorem, there exists a countable subgroup $\Gamma_0\subseteq \Gamma$ such that $f\in RI_{\Lambda_0}$, 
for $\Lambda_0=\Lambda\cap \Gamma_0$, and there exists a projection from $L^\infty_{\Lambda_0}$ onto 
${\cal C}_{\Lambda_0}$.\par
Let $H$ be the annihilator of $\Gamma_0$; that is, $H$ is the following closed subgroup of $G$:
$$H= \Gamma_0^\perp = \{ x\in G\,;\ \gamma(x)=1,\ \forall \gamma \in \Gamma_0\}.$$
The quotient group $G/H$ is metrizable since its dual group $\Gamma_0$ is countable. Let $\pi_H$ denote the 
quotient map from $G$ onto $G/H$. It is known that that the map $g\mapsto  g\circ \pi_H$ gives an isometry 
from $L^\infty_{\Lambda_0}(G/H)$ onto $L^\infty_{\Lambda_0}(G)$ sending ${\cal C}_{\Lambda_0}(G/H)$ 
onto ${\cal C}_{\Lambda_0}(G)$. \par
In order to finish our reduction to the metrizable case we only have to see that this isometry sends
$RI_{\Lambda_0}(G/H)$ onto $RI_{\Lambda_0}(G)$. It is easy to see, via the map $g\mapsto  g\circ \pi_H$, 
that having a Riemann-integrable function $g\colon G/H\to \C$ is the same as having a Riemann-integrable 
function $g\colon G\to \C$ with the property $g(x+h)=g(x)$, for every $x\in G$ and every $h\in H$. Therefore the 
above isometry sends $RI_{\Lambda_0}(G/H)$ into $RI_{\Lambda_0}(G)$. The surjectivity of this  map is a 
consequence of the following proposition:

\begin{proposition}
Let $f\colon G\to \C$ be a Riemann-integrable function such that $\hat{f}(\gamma)=0$, for every
$\gamma\in \Gamma\setminus\Gamma_0$. Then there exists a Riemann-integrable function $g\colon  G\to \C$ 
such that:\par
a) $f=g$ almost everywhere;\par
b) $g(x)=g(x+h)$, for all $x\in G$ and for all $x\in H$.
\end{proposition}

\noindent{\bf Proof.} We can and we will assume that $f$ is in fact real valued. Take an increasing sequence 
$(K_n)_n$ of compact subsets of $G$ such that, if $B=\bigcup_n K_n$, then:\par
$i)$ $f$ is continuous at every point of $B$;\par
$ii)$ $m(G\setminus B) = 0$.\par\smallskip

Using the compactness of $K_n$ and the continuity of $f$ on $B$, one can find a neighbourhood $W_n$ of $0$
such that
\begin{equation}
|f(x) - f(x+y) | \le \frac{1}{n},\qquad \hbox{for every $x\in K_n$, and every $ y\in W_n$.}\tag*{(4)}
\end{equation}

Let $(V_n)_n$ be a decreasing sequence of open symmetric neighbourhoods of $0$ such that 
$V_n + V_n \subseteq W_n$, for every $n$. For every $n$, define $f_n$ as:
$$f_n(x)= \frac{1}{ m(V_n)} \int_{V_n} f(x-y) \, dm(y)\,,\qquad x\in G.$$
$f_n$ is a continuous function since it is the convolution of $f$ and 
$$\psi_n= \frac{1}{m(V_n)}\ind_{V_n}.$$
We also have
$$\widehat{f_n}(\gamma)=\widehat{f}(\gamma)\widehat{\psi_n}(\gamma)=0\,,
\qquad\hbox{for all } \gamma\in \Gamma\setminus\Gamma_0.$$
Then the continuous function $f_n$ only depends on the classes in $G/H$; that is,
$$f_n(x)=f_n(x+h)\,, \qquad \hbox{for all $x\in G$, all $h\in H$ and all $n$.}$$

Define
$$g(x) =  \frac{1}{ 2}\Bigl(\limsup_{n\to\infty} f_n(x) + \liminf_{n\to\infty} f_n(x)  \Bigr)\,, \qquad x\in G.$$
It is clear that $g(x)=g(x+h)$, for all $x\in G$ and for all $h\in H$. Since $V_n\subseteq W_n$, we have from 
$(4)$ that $|f_n(x)-f(x)|\le 1/n$, for all $x\in K_n$. If $x\in B=\bigcup_n K_n$, then there exists
$N$ such that $x\in K_n$, for all $n\ge N$. Therefore $|f_n(x)-f(x)|\le 1/n$, for all $n\ge N$, and
$g(x)=f(x)$. So $f=g$ almost everywhere.\par

In order to finish the proof we are going to see that every point of $B$ is a point of continuity of $g$, and so $g$ is 
Riemann-integrable. Let $x$ be in $B$. Given $\varepsilon>0$, there exists $N$ such that $1/N\le \varepsilon$
and $x\in K_n$, for all $n\ge N$. We are going to prove
\begin{equation}
|g(x)-g(x+y)|\le \varepsilon,\qquad \hbox{for every $y\in V_N$.}\tag*{(5)}
\end{equation}
So $g$ will be continuous at $x$.\par

Take $n\ge N$, and $y\in V_N$. For every $z\in V_n$ we have $x+y+z\in W_N$, and $|f(x)-f(x+y+z)|\le 1/N$. 
By the definition of $f_n$ we get $|f(x)-f_n(x+y)|\le 1/N$, for every $n	\ge N$. Then we obtain $(5)$ easily, since 
$f(x)=g(x)$.\hfill\qed
\bigskip

\noindent{\bf Remarks.} 1) Actually the proof shows that if $\Lambda$ is a Lust-Piquard set and if there exists a 
surjective projection $Q\colon L^\infty_\Lambda\to {\cal C}_\Lambda$ which commutes with translations, then 
$\Lambda$ is a Rosenthal set.\par
2) Talagrand's work \cite{Ta2} uses Martin's axiom, and, in \cite{Ta3}, another axiom,called $L$. But these axioms 
do not intervene in the results we use (they are needed to obtain Riemann-integrability from the weak measurability
of translations: see \cite{Ta3}, Theorem 15-4).\par
3) F. Lust-Piquard and W. Schachermayer (\cite{LP-S}, Corollary IV.4 and Proposition IV.15; see also 
\cite{GGMS}, Theorem V.1, Corollary VI.18, and Example VIII.10) showed that if 
$L^1(G)/L^1_{\Gamma\setminus(-\Lambda)}$ does not contain $\ell_1$ (which is equivalent to 
$L^\infty_\Lambda$ having the weak Radon-Nikodym property \cite{Ta3}, Corollary (7-3-8)), then
$L^\infty_\Lambda=RI_\Lambda$. Hence $\Lambda$ must be a Rosenthal set if $L^\infty_\Lambda$ has the 
weak Radon-Nikodym property and there exists a projection from $L^\infty_\Lambda$ onto ${\cal C}_\Lambda$. 
However, a direct proof is available. For a more general result, see \cite{GGMS}, Example following 
Proposition VII.6.\par
4) The first part of the proof is the same as the one used by A. Pe{\l}czy\'nski (\cite{P}, Cor. 9.4 (a)) to show that
$A({\D})={\cal C}_\N$ is not complemented in $H^\infty=L^\infty_\N$.

\medskip
\noindent{\bf Question.} When $\Lambda$ is not a Rosenthal set, or, merely when ${\cal C}_\Lambda$ contains 
$c_0$, how big can $L^\infty_\Lambda/{\cal C}_\Lambda$ be?

\noindent{\it
P. Lef\`evre et D. Li, Universit\'e d'Artois,
Laboratoire de Math\'ematiques de Lens,
Facult\'e des Sciences Jean Perrin,
Rue Jean Souvraz, S.P.\kern 1mm 18,
62\kern 1mm 307 LENS Cedex,
FRANCE\par
\noindent pascal.lefevre@euler.univ-artois.fr\par
\noindent daniel.li@euler.univ-artois.fr
\smallskip

\noindent
H. Queff\'elec,
Universit\'e des Sciences et Techniques de Lille,\\
Laboratoire A.G.A.T.,
U.F.R. de Math\'ematiques,\\
59\kern 1mm 655 VILLENEUVE D'ASCQ Cedex,
FRANCE\par
\noindent queff@agat.univ-lille1.fr
\smallskip

\noindent
Luis Rodr{\'\i}guez-Piazza, Universidad de Sevilla, Facultad de
Matematicas, Dpto de An\'alisis Matem\'atico, Apartado de Correos 1160,
41\kern 1mm 080 SEVILLA, SPAIN\par
\noindent piazza@us.es\par}

\end{document}